\documentclass[12pt]{amsart}
\usepackage[utf8]{inputenc}
\usepackage{amsmath,amssymb,amsthm}
\usepackage{tikz,tikz-cd}
\usepackage{lipsum} 

\DeclareMathOperator{\Pic}{\mathrm{Pic}}

\newtheorem{theorem}[equation]{Theorem}
\newtheorem{question}[equation]{Question}
\newtheorem{lemma}[equation]{Lemma}
\newtheorem{corollary}[equation]{Corollary}
\theoremstyle{definition}
\newtheorem{definition}{Definition}[section]

\theoremstyle{remark}
\newtheorem{remark}[equation]{Remark}

\begin{document}

\begin{abstract}
We prove that if $(X, B+\mathbf{M})$ is a generalized klt pair with $K_X+B+\mathbf{M}_X$ nef and abundant, then $K_X+B+\mathbf{M}_X$ is semiample.
\end{abstract}

\title[Semiampleness]{Semiampleness for generalized pairs}
\author{Priyankur Chaudhuri}
\address{School of Mathematics, Tata Institute of Fundamental Research, Homi Bhabha Road, Colaba,
Mumbai 400005}
\email{pkurisibang@gmail.com}

\maketitle
\section{Introduction} A nef divisor $L$ on a normal projective variety is called \emph{abundant} if its numerical dimension (defined as the maximum power of $c_1(L)$ which is numerically nontrivial) equals its Iitaka dimension. Kawamata \cite{Ka} (see also \cite{Fu}) proved that if $(X,B)$ is a klt pair with $K_X+B$ nef and abundant, then it is semiample. From the point of view of the abundance conjecture, this is an important result. It is natural to wonder if one can prove similar results for generalized klt pairs $(X, B+\mathbf{M})$. More precisely, 

\begin{question} \label{q}
If $(X, B+ \mathbf{M})$ is a generalized klt pair with $K_X+B+\mathbf{M_X}$  nef and abundant, then is it semiample?
\end{question}

This question was investigated in \cite{PC} where we obtained a generalized canonical bundle formula for sub-klt pairs \cite[Theorem 15]{PC} and used it to answer the above question positively in case $K_X+B \geq 0$ and $\mathbf{M}_X$ is nef \cite[Corollary 20]{PC}. \\

In this article, we obtain an affirmative answer to question \ref{q}. More generally, we deduce the following global generation criterion:

\begin{theorem}
    Let $(X,B+\mathbf{M})$ be a generalized klt pair, $H$ a nef $\mathbb{Q}$-Cartier divisor on $X$ such that $H-(K_X+B+\mathbf{M}_X)$ is nef and abundant. Suppose also that $\kappa(aH-(K_X+B+\mathbf{M}_X)) \geq 0$ and 
$ \nu(aH-(K_X+B+\mathbf{M}_X)) =\nu(H-(K_X+B+\mathbf{M}_X))$
for some $a \in \mathbb{Q}_{>1}$. Then $H$ is semiample. 
\end{theorem}

We achieve this by proving a modified version (see Theorem \ref{main}) of Kawamata's characterization of nef and abundant divisors \cite[Prop 2.1]{Ka} and apply Filipazzi's generalized canonical bundle formula \cite[Theorem 1.4]{Fi} to this setting to obtain b-nefness of the moduli part. We then use \cite[Lemma 11]{PC} to conclude.

\section{Preliminaries}
\begin{definition}
[Nef reduction and nef dimension \cite{8a}] Let $X$ be a normal projective variety and let $L \in \Pic X$ be nef. Then there exists a dominant rational map $ \phi: X \dashrightarrow  Y$ with connected fibers which is proper and regular over an open subset of $Y$ (i.e. there exists $V \subset Y$ nonempty open such $\phi$ restricts to a proper morphism $ \phi|_{\phi^{-1}(V)}:\phi^{-1}(V) \rightarrow V$) where $Y$ is also normal projective such that

\begin{enumerate}
\item If $F \subset X$ is a general compact fiber of $\phi $ with $\dim F= \dim X - \dim Y$, then $L|_{F} \equiv 0 $.

\item If $x \in X$ is a very general point and $C \subset X$ a curve passing $x$ such that $\dim (\phi(C)) >0 $, then $(L \cdot C) > 0$.

\end{enumerate}

$\phi$ is called the \textbf{nef reduction map} of $L$ and $\dim Y$ the \textbf{nef dimension} $n(L)$ of $L$. 

\end{definition}

\begin{remark} Note that if $\phi: X \dashrightarrow Y$ is proper and regular over an open subset of $Y$, then there exists a resolution $\hat{\phi}:\hat{X} \rightarrow Y$ of $\phi$ such that the exceptional divisor $Ex (\hat{X}/X)$ is $\hat{\phi}$-vertical. For example, $\hat{X}$ can be chosen to be the normalization of closure of the graph of $\phi$.
\end{remark}

\begin{definition}
[Generalized pairs and their singularities(see \cite{Bi})] A \emph{generalized sub-pair} $(X,B+\mathbf{M})$ consists of a normal projective variety $X$, a $\mathbb{Q}$-divisor $B$ and a $\mathbb{Q}$-b-divisor $\mathbf{M}$ on $X$ such that:
\begin{enumerate}
\item $K_X+B+\mathbf{M}_X$ is $\mathbb{Q}$-Cartier.
\item $\mathbf{M}$ is b-nef (A b-divisor \cite{Co} is called b-nef if it descends to a nef divisor on a birational model of $X$).
\end{enumerate}
When $B \geq 0$, we drop the prefix sub.\\

Let $(X,B)$ be a generalized sub-pair and $Y \xrightarrow{\mu} X$ a higher birational model of $X$. Let $B_Y$ be defined by $K_Y+B_Y+\mathbf{M}_Y = \mu^*(K_X+B+\mathbf{M}_X)$. We say that $(X,B+\mathbf{M})$ is \emph{generalized sub-klt} (resp. \emph{generalized sub-lc}) if every coefficient of $B_Y$ is less than $1$ (resp. $\leq 1$).
\end{definition}

\section{Nef and abundant divisors}

The following two lemmas will be used in the sequel. They are both well known. However in order to avoid circular logic, I prove them avoiding the use of \cite[Prop 2.1]{Ka}.

\begin{lemma} \label{ineq} If $L$ is a nef divisor on a normal projective variety $X$, then $n(L) \geq \nu (L) \geq \kappa (L)$. 

\begin{proof}
   That $\nu(L) \geq \kappa(L)$ is well known \cite[Prop 2.7 (6)]{Nak}. Suppose $\nu(L) > n(L)$. Then $ L^{n(L)+1} \not \equiv 0$. Let $ \phi: X \dashrightarrow Y$ be the nef reduction map of $L$, so that $\dim Y = n(L)$. Let $X^{'}$ be a smooth resolution of indeterminacy of $\phi $ with $\phi^{'}:X^{'} \rightarrow Y$ the induced morphism to $Y$. By \cite[Lemma 3.1]{LP}, there exists a birational morphism $Y^{''} \rightarrow Y$ from a smooth projective variety  such that letting $X^{''}$ denote a desingularization of the main component of $X^{'} \times _Y Y^{''}$ and $\phi^{''} : X^{''} \rightarrow Y^{''} $, $\mu: X^{''} \rightarrow X$ the induced morphisms, we have $L^{''} := \mu^*L \equiv (\phi^{''})^*(M)$ for some $M \in \Pic Y$. Note that $ n(M)= \dim Y = n(L)$. But then $(L^{''})^{n(L)+1} \equiv  (\phi^{''})^*(M)^{\dim Y+1} \equiv 0$. This proves the lemma. 
\end{proof}
\end{lemma}

\begin{lemma} \label{bigness}
    Let $L$ be a nef and abundant divisor on a normal projective variety $X$ with nef dimension $n(L) = \dim X$. Then $L$ is big.

    \begin{proof}
    Suppose $\kappa(L) < \dim X$. Let 

    \begin{center}
\begin{tikzcd} 
X^{'} \arrow[d,"\mu"] \arrow[dr, "\phi"] \\
X \arrow[r, dotted] & Y

\end{tikzcd}
\end{center}

be a smooth resolution of indeterminacy of the Iitaka fibration of $L$ and let $F \subset X^{'}$ be a general fiber of $\phi$. Note that $\dim F > 0$ by assumption and $ n(L^{'}|_F) = \dim F$ by \cite[Lemma 2.10]{LP} where $L^{'} := \mu^* (L)$. Suppose $\kappa(L)= \nu(L)= k = \dim Y$. We then have 
\begin{center}
$mL^{'} = \phi^*(A)+E$
\end{center}
for some $ m \in \mathbb{N}$ where $E \geq 0$ and $A$ is ample on $Y$. We use a trick of Kawamata \cite{Ka} to show that $E$ is $\phi$-vertical. We have $ L^{'k+1} \equiv (\phi^*A)^{k+1} \equiv 0$. Then we have 

\begin{center}
    $0 \equiv(mL^{'})^{k+1}- ( \phi^*A)^{k+1} \equiv ( mL^{'} - \phi^*A)[(mL^{'})^k +(mL^{'})^{k-1}\cdot \phi^*A+ \cdots + (\phi^*A)^k] \equiv E \cdot [(mL^{'})^k +(mL^{'})^{k-1}\cdot \phi^*A+ \cdots (\phi^*A)^k]$
\end{center}
Now $L^{'}$ and $\phi^*A$ both being nef, the intersections of their powers with $E$ are all numerically effective cycles and hence are all numerically trivial. In particular, $E \cdot (\phi^*A)^k \equiv 0$ which shows that $E$ is $\phi$-vertical. Thus, $mL^{'}|_F = \phi^*A|_F \sim 0$ which contradicts the fact that $n(L^{'}|_F) = \dim F >0$. Thus $L$ is big.

    \end{proof}
\end{lemma}

The following modification of Kawamata's theorem \cite[Prop 2.1]{Ka} will be indispensible to our needs. Our proof is somewhat different from his. One can also get this by suitably modifying Kawamata's original proof.

\begin{theorem} \label{main}
    Let $X$ be a normal projective variety and $L \in \Pic X$ nef and abundant. Then there exist a birational morphism $ \mu: X^{'} \rightarrow X$ from a normal projective variety and a contraction morphism $f: X^{'} \rightarrow Y$ to a smooth projective variety such that the following conditions are satisfied:

    \begin{enumerate}
        \item The exceptional divisor $Ex (X^{'}/X)$ does not dominate $Y$.

        \item There exists a nef and big divisor $D$ on $Y$ such that $\mu^*L \sim_{\mathbb{Q}} f^*D$.
    \end{enumerate}

    \begin{proof}
    
       If $n(L)= \dim X$, then we are done by lemma \ref{bigness}. So we may assume $n(L)< \dim X $. 
       Let $ X \dashrightarrow Y_1$ be the nef reduction map of $L$. Let $X_1$ be the normalization of the closure of its graph and $f_1:X_1 \rightarrow Y_1$ the induced morphism. Then the pullback of $L$ to $X_1$ is numerically trivial on a general fiber of $f_1$. By \cite[Lemma 3.1]{LP}, there exists a birational morphism $Y \rightarrow Y_1$ from a smooth projective variety such that letting $X^{'}$ denote the normalization of the main component of $X_1 \times _{Y_1}Y$ and $ \mu: X^{'} \rightarrow X$ and $f: X^{'} \rightarrow Y$ the induced morphisms, we have $ \mu^*L \equiv f^*M$ for some nef divisor $M$ on $Y$. Note that $n(M) =\dim Y$ and $Ex (X^{'}/X)$ does not dominate $Y$.\\

       Letting $F$ denote a general fiber of $f$, note that $\kappa(\mu^*(L)|_F)\leq 0$. If $ \kappa(\mu^*(L)|_F) =-\infty $, then $F$ is contained in the stable base locus $B(\mu^*L)$ which is impossible since $F$ is a general fiber. Thus $\kappa(\mu^*(L)|_F) =0$ and hence $ \mu^*(L)|_F \sim _{\mathbb{Q}}0$. Let $U \subset Y$ be an open subset over which $f$ is flat and such that $ \mu^*(L)|_{X_y} \sim _{\mathbb{Q}} 0$ for all $y \in U$. Let 
       \begin{center}
           $K_m:= \{y \in U| h^0(m \mu^*L|_{X_y}) \geq 1\}$.
       \end{center}

       By upper semicontinuity theorem, $K_m$ is closed in $U$ for all $ m \in \mathbb{N}$. Note that $U = \bigcup _{m \in \mathbb{N}}K_m$. Now since a variety/$\mathbb{C}$ can not be a countable union of proper closed subsets, $U=K_m$ for some $m \in \mathbb{N}$. Note that since $\mu^*L$ is numerically trivial over $Y$, it follows that $m\mu^*(L)|_{X_y} \sim 0$ for all $y \in U$. So there exists an $f$-vertical divisor $E$ on $X$ such that $m\mu^*(L)-E \sim 0$. Since $f(E)$ is a proper subvariety of $Y$, there exists a divisor $D$ on $Y$ such that $E-f^*D$ supports no fibers over codimension $1$ points of $Y$. By the generalized negativity lemma \cite[Lemma 3.3]{Bi2}, we have $E-f^*D \geq 0$ and $-E+f^*D \geq 0$. Thus $E=f^*D$. We thus conclude that $\mu^*L \sim_{\mathbb{Q}}f^*D$. Now since $D$ is nef and abundant of maximal nef dimension, we conclude that $D$ is big by Lemma \ref{bigness}.

        \end{proof} 
\end{theorem}

       \begin{corollary}\label{main}
            Let $(X,B+\mathbf{M})$ be a generalized klt pair, $H$ a nef $\mathbb{Q}$-Cartier divisor on $X$ such that $H-(K_X+B+\mathbf{M}_X)$ is nef and abundant. Suppose also that $\kappa(aH-(K_X+B+\mathbf{M}_X)) \geq 0$ and 
$ \nu(aH-(K_X+B+\mathbf{M}_X)) =\nu(H-(K_X+B+\mathbf{M}_X))$
for some $a \in \mathbb{Q}_{>1}$. Then $H$ is semiample. 

\begin{proof}
    Consider the diagram $X \xleftarrow {\mu} X^{'} \xrightarrow {f} Y$ associated to $N:=H-(K_X+B+\mathbf{M}_X)$ such that $N^{'} :=\mu^*N \sim_{\mathbb{Q}} f^*D$ as given by the previous theorem.  Letting $H^{'}:= \mu^*H$, we have $[(a-1)H^{'}+N^{'}]^{k+1}\equiv 0$ where $k= \kappa (N)$. Both the divisors being nef, we have $H^{'} \cdot N^{'k} \equiv 0 $. From now on, we may replace $(a-1)H^{'}$  by $H^{'}$.\\
    
    Writing $mD= A+E$ for some $m \in \mathbb{N}$ where $A$ is very ample and $E \geq 0$, arguing as in the proof of \cite[Prop 2.61]{KM}, we can show that $m^kD^k \geq A^k$. Therefore, $N^{'k} \geq \frac{1}{m^k} (f^*A)^k$ (as cycles). Since $H^{'} \cdot N^{'k} \equiv 0$, we conclude that $(f^*A)^k \cdot H^{'} \equiv 0$. Thus if $F$ denotes a general fiber of $f$, $H^{'}|_F \equiv 0$. We also have $N^{'}|_F \equiv 0$. Since $\kappa(H^{'}+N^{'}) \geq 0$, $(H^{'}+N^{'})|_F \sim _{\mathbb{Q}} 0$. Arguing as in the proof of the above theorem (possibly after a birational base change of $f$ which preserves the $f$-verticality of the $\mu$- exceptional divisor), there exists a divisor $D^{'}$ on $Y$ such that $H^{'}+N^{'} \sim _{\mathbb{Q}} f^*(D^{'})$. Thus $\mu^*(K_X+B+\mathbf{M}_X) \sim_{\mathbb{Q},f}0$. Consider the induced generalized canonical bundle formula

    \begin{center}
        $(X^{'},B_{X^{'}}+\mathbf{M}_{X^{'}}) \xrightarrow{f} Y $ 
        
        \end{center}
        where $K_{X^{'}}+B_{X^{'}}+\mathbf{M}_{X^{'}}:= \mu^*(K_X+B+\mathbf{M}_X)$ and $B_{X^{'}}$ is defined by the requirement $\mu_* \mathbf{M}_{X^{'}}= \mathbf{M}_X$. Now $B_{X^{'}}$ differs from $\mu_*^{-1}B$ by some $\mu$-exceptional divisor. The exceptional divisor $Ex(X^{'}/X)$ is $f$-vertical. Therefore, $(X^{'}, B_{X^{'}}+\mathbf{M}_{X^{'}})$ is generalized klt over the generic point of $Y$. Then by \cite[Theorem 1.4]{Fi}, we can write $K_{X^{'}}+B_{X^{'}}+\mathbf{M}_{X^{'}} \sim _{\mathbb{Q}} f^*(K_Y+B_Y+\mathbf{M}_Y)$ where $B_Y$ is the discriminant and $\mathbf{M}_Y$ is the moduli b-divisor of the generalized klt-trivial fibration $f$ such that $\mathbf{M}_Y$ is b-nef. Recall that $B_Y$ is defined in the following way: \\
        
        For any prime divisor $P$ on $Y$, let $a_P :=sup \{t|(X^{'},B_{X^{'}}+tf^*P+\mathbf{M}_{X^{'}})$ is generalized sub-lc over the generic point of $P$\} and define $B_{Y^{'}}:= \Sigma_P(1-a_P)P$. Then $H$ is semiample by \cite[Lemma 11]{PC}.
    
\end{proof}
       \end{corollary}

       The following consequence of Corollary \ref{main} is the analogue of \cite[Corollary 3.5]{CKP} for generalized pairs.
       \begin{corollary}\label{num}
           Let $(X,B+\mathbf{M})$ be a generalized klt pair. Then $K_X+B+\mathbf{M}_X$ is numerically equivalent to a nef and abundant divisor if and only if it is numerically equivalent to a semiample divisor.

           \begin{proof}
               One direction is obvious. Let $K_X+B+\mathbf{M}_X$ be numerically equivalent to a nef and abundant divisor. Then there exists a numerically trivial divisor $N$ on $X$ such that $K_X+B+\mathbf{M}_X+N$ is nef and abundant. We apply Corollary \ref{main} to the generalized klt pair$(X,B+\mathbf{M+N})$ where the trace of $\mathbf{M+N}$ on a higher model $X^{'} \xrightarrow{\pi}X$ is $\mathbf{M}_{X^{'}}+\pi^*N$. We conclude that$K_X+B+\mathbf{M}_X+N$ is semiample and thus $K_X+B+\mathbf{M}_X$ is numerically equivalent to a semiample divisor.
           \end{proof}
       \end{corollary}
       \begin{remark}
          Corollary \ref{num} is the analogue of Kawamata's theorem \cite[Theorem 1.1]{Ka} in the setting of generalized abundance-type conjectures(see \cite[page 2]{LP} \cite[Conjecture 1.3]{HL}).
       \end{remark}

       \section{Acknowledgements}

       I have greatly benefited from several conversations with Omprokash Das on the contents of this note. I also thank Christopher Hacon for an email conversation and the Tata Institute of Fundamental Research for providing excellent working conditions. Partial financial support was provided by the Department of Atomic Energy, India under project no. 12-R\&D-TFR-5.01-0500 during this work.

\end{document}